\documentclass[12pt]{article}
 
\usepackage{graphicx}

\usepackage{amsfonts}
\usepackage{amssymb}
\usepackage{amsthm}
\usepackage{amsmath}
\usepackage{enumerate}
\usepackage{rotating}
\usepackage{graphicx}
\usepackage{bm}
\usepackage{graphicx}
\usepackage{color}
\newcommand{\beqn}{\vspace{-0.25cm}\begin{eqnarray*}}

\newcommand{\eeqn}{\end{eqnarray*}}
\newcommand{\bneqn}{\vspace{-0.25cm}\begin{eqnarray}}
\newcommand{\eneqn}{\end{eqnarray}}

\newcommand{\parens}[1]{\left(#1\right)}

\newcommand{\Prob}{\mathbb{P}}
\newcommand{\Exp}{\mathbb{E}}
\newcommand{\Expv}{\mathcal{E}}

\newcommand{\bfZ}{\mathbf{Z}}
\newcommand{\bfQ}{\mathbf{Q}}
\newcommand{\bfG}{\mathbf{G}}
\newcommand{\bfM}{\mathbf{M}}
\newcommand{\bfI}{\mathbf{I}}
\newcommand{\bfLambda}{\mathbf{\Lambda}}

\newcommand{\bfDelta}{\mathbf{G}^{-1}}
\newcommand{\bfT}{\mathbf{T}}
\newcommand{\bfW}{\mathbf{W}}
\newcommand{\bfV}{\mathbf{V}}
\newcommand{\bfe}{\mathbf{e}}

\newcommand{\bfx}{\mathbf{x}}

\newcommand{\bfw}{\mathbf{w}}

\newcommand{\bftau}{\mathbf{\tau}}
\newcommand{\bfzero}{\mathbf{0}}

\newcommand{\Fb}{{\overline F}}
\newcommand{\dd}{{\mathrm{d}}}
\newcommand{\e}{{\mathrm{e}}}

\newtheorem{theorem}{Theorem}
\newtheorem{proposition}[theorem]{Proposition}
\newtheorem{lemma}[theorem]{Lemma}
\newtheorem{corollary}[theorem]{Corollary}
\newtheorem{example}[theorem]{Example}
\newtheorem{remark}[theorem]{Remark}

\def\tcr{\textcolor{black}}\def\tcrr{\textcolor{black}}
\def\tcb{\textcolor{black}}

\title{Stability and busy periods in a multiclass queue with
state-dependent arrival rates
}
\author{ Philip A. Ernst, S\o{}ren Asmussen, and
        John J. Hasenbein}
\begin{document}
\maketitle

\maketitle

\begin{abstract}
We introduce a multiclass single-server queueing system in which the arrival
rates depend on the current job in service. The system is characterized by a matrix of arrival rates in lieu of a vector of arrival rates. Our proposed model departs from existing state-dependent queueing models in which the parameters depend primarily on the number of jobs in the system rather than on the job in service.  We formulate the queueing model and its corresponding fluid model and proceed to obtain the necessary and sufficient conditions for stability via fluid models. Utilizing the natural connection with the multitype Galton-Watson processes, 
the Laplace-Stieltjes transform of busy periods in the system is given. We conclude with 
tail asymptotics for the busy period 
for heavy-tailed service time distributions for the regularly varying case. 
\end{abstract}
Keywords: Busy periods; fluid models; multiclass queues; regular variation; stability; state-dependent arrival rates.

\section{Introduction}
\setcounter{equation}{0}
We introduce a multiclass single-server queueing system in which the arrival
rates depend on the current job in service. The system is characterized by a matrix of arrival rates instead of a vector of arrival rates. The proposed model departs from existing state-dependent models in the literature in which the parameters depend primarily on the number of jobs in the system (see Bekker et al.~\cite{Bekker}, Cruz and Smith~\cite{Cruz}, Jain and Smith~\cite{Jain}, Perry et al.~\cite{Perry}, and Yuhaski and Smith~\cite{Yuh}, among other sources) rather than the job in service. \\
\indent Our model is motivated by two practical queueing considerations. The first
is a multiclass queueing system in which the arriving customer can observe only the class of the customer in service and no other characteristics of the queue. This information informs the customer's decision to either join or leave the queue. The second concerns local area networks with a central server in which $K$ clients generate requests
at individual Poisson rates $\mu_i$. Often, a client does not generate requests when
a previous request is being handled by the server. Further, it is conceivable that groups of clients working together may influence each other's Poisson rate. To the best of our knowledge, this simple yet potentially very useful queueing model has never appeared in the literature. This serves as our primary motivation for the manuscript.\\
\indent The remainder of the work is structured as follows.  We formulate the queueing model in Section \ref{sec2} and its corresponding fluid model in Section \ref{sec3}.  In Section \ref{sec4}, we obtain the necessary and sufficient conditions for stability via fluid models.  Through the natural connection with the multitype Galton-Watson processes, we characterize
the Laplace-Stieltjes transform of busy periods in the system in Sections \ref{sec5} and \ref{sec6}. Section \ref{sec7} concerns tail asymptotics on the busy period in the case of heavy-tailed service time distributions. Section \ref{sec8} offers a brief conclusion and presents ideas for future work.

\section{The queueing model} \label{sec2}

Consider a multiclass single-server queue with $K$ classes of jobs, each arriving according to independent counting processes. \tcr{We assume that only one job may be serviced at a time}. Let the arrival rate depend \textit{on the class of the job in service}.  Let the arrival rate depend \textit{on the class of the job in service}. If the server is serving a job of class $i$, the arrival rate of class $j$ jobs is $\lambda_{ij}$, $i,j=1,\dots,K$.
The matrix of arrival rates is defined as $\bfLambda = (\lambda_{ij})$, $i,j = 1, \ldots, K$. 
If there is no job in service, then the arrival rate of class $j$ jobs is defined as $\lambda_{0j}$, 
$j=1,\dots,K$.
The arrival mechanism is described more precisely with dynamical equations in Section \ref{sec3}.

We proceed to set notation. Let 
$\bar{\lambda}^i = \sum_{j=1}^K \lambda_{ij}$ 
for each $i=1,\dots,K$.
Service times for class $i$ jobs are assumed to be i.i.d.\ with distribution function $F_i$, 
$i=1,\dots, K$.  Let $S_i$ be a generic service time for class $i$ jobs, with  $\Exp[S_i] = m_i =\mu_i^{-1}$,  $i=1,\dots, K$ and $\bfG=\text{diag}(\mu_1,\mu_2,\ldots,\mu_K)$. We define the ``mean offspring matrix'' 
to be $\bfM = \bfG^{-1} \bfLambda $ (here, the $ij$th element $\lambda_{ij}m_i$
is the mean number of arriving class $j$ customers during service of a class $i$ customer).
By definition, all the elements of $\bfM$ are non-negative, and this is enough to ensure that 
the dominant eigenvalue $\rho(\bfM)$ is real and positive, cf.~\cite{Gant}.
For some results,  more restrictive conditions on $\bfM$ will be required. Further, let $\psi_i$ denote the Laplace-Stieltjes transform (LST) of  $S_i$, $i=1,\dots, K$, respectively,  that is,  $\psi_i(s) = \Exp[\e^{-s S_i}] = 
\int_0^{\infty} \e^{-st} \dd F_i(s)$ for $s>0$. We let $Q_i$ denote the steady-state number of class $i$ jobs in the system,  $i=1,\dots, K$, and let $\bfQ= (Q_1,\dots, Q_K)$.
Each state of the system takes nonnegative integer values, that is, $\bfx = (x_1,\dots,x_K) \in \mathbb{Z}^K_+$.

The service disciplines we consider are non-idling, i.e., jobs must be served using the full
capacity of the server whenever there are jobs in the system. 
\tcr{Our results on stability and the busy period are independent of the particular
(non-idling) scheduling policy employed in the system.}


\section{Queueing and fluid dynamics} \label{sec3}\setcounter{equation}{0}
\subsection{Queueing dynamical equations} 

We now precisely define the arrival mechanism. 
For $i \in \{1, \ldots, K\}$ and $t \ge 0$, $Q_i(t)$ denotes the number
of class $i$ jobs in the system at time $t$, whether in service or in the queue. 
Similarly, let $T_i(t)$ denote
the amount of time that has been devoted to serving class $i$ jobs
in $[0,t]$. Further, let $A_i(t)$ and $D_i(t)$ be, respectively,
the total number of class $i$ jobs that have arrived and departed from the
system in $[0,t]$. We then have the following input-output equation for
each class $i$ job
\begin{equation}
Q_i(t) = Q_i(0) + A_i(t) - D_i(t).
\end{equation}

For each class $i$, the counting process $\Expv^j_i(t)$ is the number of class $i$
jobs that arrive during the first $t$ time
units devoted to processing class $j$. $\Expv^0_i(t)$ counts the number of class $i$
arrivals during the first $t$ time units that no job is being processed at the server. 
The total number of class $i$ arrivals in $[0,t]$ is then given
by
\begin{equation} 
A_i(t) = \Expv^0_i\bigl(T_0(t)\bigr) +   \sum_{j=1}^N \Expv^j_i \bigl(T_j(t)\bigr),
\end{equation}
where the counting processes $\Expv^j_i$ for \tcr{$i=1,\dots, K,$ and $j = 0, \dots, K$} are assumed to be mutually independent. 

As for the service processes, for each $i$, $1 \le i \le K$ and positive integer
$n$, we let $V_i(n)$ denote the total service requirement for the first $n$
class $i$ jobs. Assuming an HL service discipline, we have that
\begin{equation} 
V_i(D_i(t)) \le T_i(t) \le V_i(D_i(t)+1) 
\end{equation}
for each $t \ge 0$ and $1 \le i \le N$.

We define the workload in the system at time $t$ to be 
\begin{equation} 
W(t) = \sum_{i=1}^K V_i(A_i(t) + Q_i(0)) - \sum_{i=1}^K T_i(t),
\end{equation}
and the cumulative idle time process to be
\begin{equation} 
Y(t) = t - \sum_{i=1}^K T_i(t).
\end{equation} 
It is important to note that $Y$ is a non-decreasing function. We assume that the queueing
policy is non-idling, which specifically means that $Y$ can increase only
when $W(t)=0$. More precisely,
$$ \int_0^{\infty} W(t) \, dY(t) = 0.$$

\subsection{Fluid model}
For purposes of determining the stability conditions of a more general 
version of our model, we formulate a fluid network version of the model. For references to important definitions and results in the fluid model literature, we refer the reader to Bramson~\cite{Bramson} and Gamarnik~\cite{Gamarnik}.

For $i \in \{1, \ldots, K\}$ and $t \ge 0$, let $Q_i(t)$ denote the amount
of fluid of class $i$ in the system at time $t$. Similarly, let $T_i(t)$ denote
the amount of time that has been devoted to serving class $i$ fluid
in $[0,t]$. We also define $A_i(t)$ and $D_i(t)$ which are, respectively,
the total amount of class $i$ fluid that has arrived and departed from the
system in $[0,t]$. We then have the following standard
equation:
\begin{equation}
Q_i(t) = Q_i(0) + A_i(t) - D_i(t),
\end{equation}
for each $i \in \{1, \ldots, K\}$ and $t \ge 0$.
The departure processes in this system also obey the standard relation
$D_i(t) = \mu_i T_i(t)$ for all $t \ge 0$. 

The unusual feature of our model lies in the arrival process, which
is dependent on the current class in service. In the queueing model,
processor sharing is not allowed. Hence, there is (at most) one class
in service at any given time and ``the customer in service'' is defined 
unambiguously. Here, we provide a more general formulation that
reduces to the queueing model presented in earlier sections, under appropriate restrictions on the 
allowable queueing disciplines. First, we recall the usual condition
\begin{equation}
\sum_{i=1}^N \dot{T}_i(t) \le 1,
\end{equation}
which simply indicates that the server cannot devote more than 100\%
of its time to serving fluids of all classes. Since we assume that
the queueing discipline is non-idling,  $\sum_{i=1}^N \dot{T}_i(t) = 1$
whenever there is a positive amount of fluid in the system. We also define
the idle time in $[0,t]$ to be:
$$ Y(t) = t - \sum_{i=1}^N T_i(t).$$

Note that the current arrival rate of class $j$ fluid is given by
$\dot{A}_j(t)$. In the queueing model, if a job of class $i$ is in service
then the arrival rate of class $j$ jobs is $\lambda_{ij}$. Let 
$\mathbf{\lambda}_j$ be the column vector $(\lambda_{1j}, \ldots
\lambda_{Nj})^\bot$ and let $\dot{\mathbf{T}}(t)$ be the column vector
$\bigl(\dot{T}_1(t), \ldots, \dot{T}_N(t)\bigr)^\bot$, where ${}^\bot$ means transposition. We define the fluid arrival
rate of class $j$ to be
\begin{equation} 
\label{arrivalrep}
\dot{A}_j(t) = \lambda_{0j} \dot{Y}(t) + \lambda^\bot_j \dot{\mathbf{T}}(t).
\end{equation}
In particular, when there is fluid in the system, the class $j$ arrival
rate is a convex combination of the elements of $\mathbf{\lambda}_j$.
If we restrict to policies in which only one class can be served at
any time, then equation (\ref{arrivalrep}) assigns an arrival rate of $\lambda_{ij}$
to class $j$ fluid when class $i$ fluid is in service. Note that this concurs with
the queueing model formulation. Combining the above, we have
\bneqn \label{eq9}
Q_j(t) & = & Q_j(0) + \int_0^t  (\lambda_{0j} \dot{Y}(u) + \lambda^\bot_j  \dot{\mathbf{T}}(u)) \; du
- \mu_j T_j(t)  \\
& = & Q_j(0) + \lambda_{0j} {Y}(t) + \lambda^\bot_j \mathbf{T}(t) - \mu_j T_j(t). \label{eq10}
\eneqn
Writing equations (\ref{eq9}) and (\ref{eq10}) in matrix form yields
\begin{equation} \label{dynamics}
\bfQ(t) = \bfQ(0)+ (\mathbf{M}^\bot-\mathbf{I})\mathbf{D}(t) +Y(t) \mathbf{\lambda_0}.
\end{equation}
We define the vector of fluid work in the system at time $t$ to be 
\begin{equation}  \label{work}
\bfW(t) = \bfG^{-1}\bfQ(t).
\end{equation}







\subsubsection{Fluid Limits}

Thus far we have described a fluid model but it remains to show
that the fluid limits of the queueing model satisfy the fluid model equations.
In this subsection only, we use a bar to denote a fluid limit.
As usual, we define the fluid
limit of the queue-length processes to be 
$$ \bar{Q}_i(t) = \lim_{n \to \infty} \frac{Q_i(nt)}{n},$$
with other fluid limits defined in an analogous manner. 
We make the usual assumptions on the stochastic primitives and initial conditions, i.e.,
for all $i$ and $j$
\begin{eqnarray}
 \lim_{n \to \infty} \frac{\Expv^j_i(nt)}{n} & = & \lambda_{ji}t   \label{fslln1}  \\
  \lim_{n \to \infty} \frac{\Expv^0_i(nt)}{n} & = & \lambda_{0i}t \label{fslln2}  \\
  \lim_{n \to \infty} \frac{V_i(n)}{n} & = & m_i \\
    \lim_{n \to \infty} \frac{Q_i(0)}{n} & = & \bar{Q}_i(0), 
\end{eqnarray}
where the convergence is almost surely, uniformly on compact sets. 
Under these assumptions, the fluid model equations
can be derived in a straightforward way from the queueing dynamical equations, since
all but the arrival rate process is identical to the standard multiclass queueing network
model. For the arrival process, we have 
\begin{eqnarray*}
\bar{A}_i(t) & = & \lim_{n \to \infty} \frac{A_i(nt)}{n} 
\ = \ \lim_{n \to \infty} \frac{\Expv^0_i(T_0(nt))}{n} + \lim_{n \to \infty}   \sum_{j=1}^N \frac{\Expv^j_i (T_j(nt))}{n} \\
& = & \lambda_{0i} \bar{Y}(t) + \lambda^\bot_i \bar{\mathbf{T}}(t).
\end{eqnarray*}

The last equality follows from assumptions (\ref{fslln1}) and (\ref{fslln2}) and
similar arguments as found in Proposition 4.12 in \cite{Bramson}.
Finally, the connection between fluid stability and
queueing network stability follow from straightforward modification of existing stability 
results, under the usual assumptions that the interarrival times for all job classes are unbounded 
and spread out. We refer the reader to Chapter 4 of Bramson \cite{Bramson} for full details. 

\section{Stability results for fluid model} \label{sec4}\setcounter{equation}{0}

In this section, we prove a number of results regarding the stability, or instability,
of the fluid model. The proofs rely on the following two observations.
\begin{enumerate}
\item $T_i(\cdot)$ is Lipschitz continuous for each $i$ and hence so is
any linear function $f$ of $(T_1, \ldots, T_K)$. Thus, $f$ is 
absolutely continuous and its derivative exists almost everywhere.
\item If $\dot{f}(t)$ exists for $t >0$, 
$t$ is called a regular point. 
\end{enumerate}
We
define $\bfe = (1, \ldots, 1)^\bot$ and assume this column vector is of size $K$. 
Finally, we set $\bf{H} = \bf G\bfM\bfG^{-1}.$ 

\begin{theorem} 
If $\rho(\textbf{M}) <1$, then 
$f(t) = \bfe^\bot (\textbf{I}-\textbf{H}^\bot)^{-1} \textbf{G}^{-1}Q(t)$
is a Lyapunov function for the fluid model. 
\end{theorem}

\begin{proof}
Note that $\rho(\textbf{H})= \rho(\textbf{M}) <1$. 
Hence, $\textbf{I}-\textbf{H}$ is an $\cal{M}$-matrix. Therefore $\textbf{I}-\textbf{H}$ is invertible
with a non-negative inverse. 

Let us assume that the fluid system starts from a non-empty state, i.e., $\bfQ(0) \not= \bfzero$. 
By the continuity of $\bfQ$, $\bfQ(t)\not= 0$ for all $t$ in some interval $[0,s)$. 
Then we have $Y(t) =0$ for all $t \in [0, s)$. Using equations (\ref{dynamics}) and (\ref{work}) we have 
$$\bfW(t) = \bfW(0) - (\textbf{I}-\textbf{H}^\bot)\bfT(t),$$
for $t \in [0, s)$.
Multiplying by $\bfe^\bot (\textbf{I}-\textbf{H}^\bot)^{-1}$ yields
$$ \bfe^\bot (\textbf{I}-\textbf{H}^\bot)^{-1}\bfW(t) = \bfe^\bot (\textbf{I}-\textbf{H}^\bot)^{-1} \bfW(0) - \bfe^\bot \bfT(t).$$
As in the statement of the theorem, set
$$ f(t) =  \bfe^\bot(\textbf{I}-\textbf{H}^\bot)^{-1} \textbf{G}^{-1}\bfQ(t)$$ 
and note that $f(t)=0$ if and only if $\bfQ(t) =\bfzero$. 
Then we have:
\begin{eqnarray}
f(t) & = & \bfe^\bot (\textbf{I}-\textbf{H}^\bot)^{-1} \bfW(t) \\
& = & \bfe^\bot (\textbf{I}-\textbf{H}^\bot)^{-1} \bfW(0) - \bfe^\bot \bfT(t).
\end{eqnarray}
 Taking derivatives, we obtain
 $$ \dot{f}(t) = -\bfe^\bot \dot{\bfT}(t) = -1,$$
 for any $t \in [0,s)$ and regular point $t$.  
 Therefore, the draining time of the system under
 any feasible policy is $$f(0) = \bfe^\bot (\textbf{I}-\textbf{H}^\bot)^{-1} \bfW(0),$$
which can be interpreted as the initial unfinished
``potential'' work, defined as the work due to the current workload
and work generated in the future by the initial
workload's ``offspring.'' The above argument implies that $\dot{f}(t)=-1$ whenever
$\bfW(t) \not= \bfzero$ and thus the system stays drained once $\bfQ(t) =\bfzero$. This completes the proof. 
 \end{proof}
The corollary below now immediately follows.
\begin{corollary} The fluid model is globally stable if $\rho(\textbf{M}) <1$. 
\end{corollary}

\subsection{Weak instability}
Next we show that the fluid model is weakly unstable  if $\rho(\textbf{M})>1$. We begin by introducing the following lemma.
\begin{lemma} \label{hithere}
Suppose $\rho(\textbf{M})=\rho(\textbf{H})>1$ and that each row of $\bfM$ has at least one strictly positive element. Then for all nonnegative
vectors $\bfT(t)>0$, $\bfV(t)=(\textbf{I}-\textbf{H})\bfT(t)$ must have some component $V_{i}(t)<0$
for some $i\in\{1,...,K\}$.
\end{lemma}

\begin{proof}
We argue to the contrary. Note that each row of $\textbf{H}$ has at least one strictly positive element, by the same assumption
on $\textbf{M}$. Also, for some $\alpha \in (0,1)$,
$\rho(\alpha \textbf{H}) = 1.$
For the sake of contradiction, assume that there exists a nonnegative
vector $\bfT(t)>0$ s.t.\ $\bfV(t)=(\textbf{I}-\textbf{H})\bfT(t)\geq0$. 
Further, define $\bfV^{'}(t) = (\textbf{I}- \alpha \textbf{H})\bfT(t)$.
We now consider:
\beqn
\bfV(t)-\bfV^{'}(t)=(\textbf{I}-\textbf{H})\bfT(t)-(\textbf{I}- \alpha \textbf{H})\bfT(t)=(\alpha \textbf{H}-\textbf{H})\bfT(t)<0.
\eeqn
The above equations imply that $\bfV^{'}(t)>\bfV(t)$ and that there exists some
$\bfT(t)>0$ with $(\textbf{I}- \alpha \textbf{H})\bfT(t)>0$. Thus $(\textbf{I}-\alpha \textbf{H})$ is semipositive, and by condition $I_{27}$ in Chapter 6 of Berman~\cite{Berman}, $(\textbf{I}-\alpha \textbf{H})$ is a non-singular $\cal{M}$-matrix. This implies that $\rho(\alpha \textbf{H}) < 1,$ yielding a contradiction. 
\end{proof}

We are now ready to prove Theorem \ref{thm2},  the main result of this subsection. 
\begin{theorem} \label{thm2}
The fluid model is weakly unstable if $\rho(\bfM)>1$ and each row of $\bfM$ has at least one strictly positive element.
\end{theorem}
\begin{proof}
Assume $\bfQ(0)=\bfW(0)=\bfzero$. Then for any $t >0$ we have by (\ref{dynamics}) and (\ref{work}) that
\beqn
\bfW(t) \ge \bfW(0) - (\textbf{I}-\textbf{H}^\bot)\bfT(t)=(\textbf{H}^\bot-\textbf{I})\bfT(t).
\eeqn
By Lemma \ref{hithere}, there exists some component of $\bfW(t)$ s.t.\ $W_{i}(t)>0$. This implies that $\bfQ(t)\neq \bfzero$ for all $t >0$. Thus the fluid model is weakly unstable. 
\end{proof}

\subsection{Weak Stability}

\begin{theorem} \label{implemma}
Suppose that $\bfM$ is an irreducible non-negative matrix.
\tcr{Then the} fluid model is weakly stable if $\rho(\bfM) \le 1$. 
\end{theorem}
\begin{proof}
It suffices to show the result for the case $\rho(\bold{M})=1$, since we have already shown that the fluid model
is ``strongly'' stable when $\rho(\bold{M})<1$.

Let $\bfQ(0)=0$. We argue to the contrary. For the sake of contradiction, let us assume that $\bfQ(t) \not= \bfzero$ for some $t >0$. Then, since $\bfQ$ is continuous,
there must be an interval $(t_1, t_2)$ with $t_2 > t_1$,  for which $\|\bfQ(t)\| > 0$ for all $t \in (t_1, t_2)$, 
$\bfQ(t_1) =0$ and $\|\bfQ(t_2) \| >0$.  In particular, we may set $t_1 = \text{inf}\{t: \bfQ(t) \not=\bfzero \}$. Now, recall that 
\begin{equation} \label{dynam1}
\bfQ(t)=(\bfM^\bot-\bfI)\bold{D}(t)+Y(t)\mathbf{\lambda_0}.
\end{equation} Since $\bold{M}$  is a positive matrix, it follows by the Perron-Frobenius Theorem that there exists a positive left (row) eigenvector $\bold{w}$  of $\bold{M}$  with $\bold{w}\bfM=\bold{w}$, $w_i >0 $  for $i \in \{1,...,K\}$. Multiplying both sides of (\ref{dynam1}) by $\bold{w}$ we obtain
$$ \bold{w}\bfQ(t)=\bold{w}\bigl[(\bfM^\bot-\bfI)\bold{D}(t)+Y(t)\mathbf{\lambda_0}\bigr]\ =\ \bold{w}Y(t)\mathbf{\lambda_0},$$
for all $t \ge 0$. Recalling $\bfQ(t_1)=\bfzero$ and $\|\bfQ(t_2)\| >0$ we have
$$ \bold{w}\bigl(Y(t_2)-Y(t_1)\bigr) \mathbf{\lambda_0} \ =\ \bold{w}\bigl(\bfQ(t_2)-\bfQ(t_1)\bigr) > 0.$$
This implies $Y(t_2) > Y(t_1)$ and thus there is positive idle time in $(t_1,t_2)$. However,
since the fluid level is positive in this entire interval, this violates the non-idling condition.
Thus such a fluid solution is not feasible. A contradiction has been reached. This concludes the proof.

\end{proof}

\section{Branching process connection} \label{sec5}\setcounter{equation}{0}

In the remainder of the paper, we investigate a special case of the multiclass
model discussed so far. In particular, we now assume that arrivals to each class from a Poisson
process, i.e., the model is an $M/G/1$ multiclass queue, rather than a $GI/G/1$ queue. 
Although more general stability conditions for the $GI/G/1$ case were proven in Section \ref{sec4}, we
begin by reproving them in the Poisson setting, by making a connection to branching
processes. There are two reasons to do this. First, the stability results arise in a somewhat more
intuitive manner using this methodology. Secondly, we find the connection to branching processes
illuminating and useful in later sections.
 
A classical tool for the simple $M/G/1$ queue and related systems is to interpret customers as individuals in a branching process,
such that the children of a customer is the number of customers arriving during his or her service. This is useful  because
the stability condition for the queueing system is the same as the condition for almost sure extinction. Carrying out the same idea
for our multiclass systems leads to a $K$-type Crump-Mode-Jagers branching process
$\{ \bfZ_n=\bigl(Z^{(1)}_n,\dots Z^{(K)}_n\bigr): n \ge 1\}$, such that the lifetime of an individual
of type $j$ has the same distribution as $S_j$. \tcr{In the results below, we consider a branching process with a single ancestor of
type $i$. Whenever $\Exp_i$ and $\Prob_i$ are used, they are with reference to the probability
measure induced by such a single ancestor.} The offspring mechanism is then
described by the probabilities
\begin{equation}
p_{ij}(k) = \Prob_i\big(Z_{1}^{(j)}= k\big) = \Prob\bigl({\text{Pois}(\lambda_{ij}S_j) = k}\bigr) = \int_0^{\infty}  \frac{(\lambda_{ij}s)^k 
\e^{-\lambda_{ij}s}}{k!} \dd F_j(s)\,.
\end{equation}
The offspring matrix $\bfM = (M_{ij})_{i,j= 1,\dots, K}$
is given by $M_{ij}=\Exp_i\bigl[Z_1^{(j)}\bigr]=\lambda_{ij}/\mu_{i}$ and is assumed
irreducible. Thus Perron-Frobenius theory applies to $\bfM$
and as before, $\rho=\rho(\bfM)$ the largest eigenvalue. Note that the $ij$th element of the matrix $\sum_{n=0}^\infty \bfM^n$ gives the expected number of
type $j$ progeny of an individual of type $i$; of course, when $\rho<1$, we have $\sum_{n=0}^\infty \bfM^n=(\bfI-\bfM)^{-1}$. 

 \subsection{Stability Conditions}
Let $|\bfZ_n|=\sum_{j=1}^KZ_n^{(j)}$ denote the total number of individuals in the $n$th generation
and $T^*$ 
the extinction time.  
Let \tcr{$\Prob_i(T^*<\infty)$ be the extinction probability of type $i$ of the branching process.}
Then, by classical results, we have the following theorem:\\

\begin{theorem}\label{thmHA}
\bneqn \label{Harris}
\Prob_i(T^*<\infty) =1,\,\, i=1,\ldots,K,\, \text{if and only if} \,\, \rho \leq 1.
\eneqn
\end{theorem}
\begin{proof}
By the classical result for the extinction time of branching processes~\cite[Chap II. Theorem 7.1]{Harris63}, if and only if $\rho \leq 1$, the total number of generations for each type is finite with probability $1$ and thus $\sum |\bfZ_n| < \infty$, which further implies $\Prob_i(T^*<\infty) = 1 $  for every $i$.

\end{proof}
\noindent Consider $K=2$. Straightforward algebra gives that $\rho \leq 1$ is equivalent to

\bneqn
\frac{\frac{\lambda_{11}}{\mu_1}+\frac{\lambda_{22}}{\mu_2}+\sqrt{\parens{\frac{\lambda_{11}}{\mu_1}-\frac{\lambda_{22}}{\mu_2}}^2+\frac{4\lambda_{12}\lambda_{21}}{\mu_1\mu_2}}}{2} \leq 1.
\eneqn

\begin{theorem}\label{thmHA2}  $\Exp_iT^*<\infty$ for all $i$ if and only if $\rho < 1$.
\end{theorem}
\begin{proof}  For a simple proof of sufficiency, assume $\rho<1$ and let 
$S_j(m;n)$ denote the lifetime of the $m$th individual of type $j$ in the $n$th generation, $\underline\mu=\min_1^K\mu_j$.
Then
\begin{align*}\Exp_iT^*\ &= \Exp_i\sum_{n=0}^\infty\sum_{j=1}^K \sum_{m=1}^{Z_n^{(j)}}S_j(m;n)\ =\ 
 \Exp_i\sum_{n=0}^\infty\sum_{j=1}^K\frac{Z_n^{(j)}}{\mu_j}\\ &\le\ 
\tcr{ {\underline\mu}^{-1} \Exp_i\sum_{n=0}^\infty\sum_{j=1}^KZ_n^{(j)}\ =\ {\underline\mu}^{-1} \sum_{n=0}^\infty\sum_{j=1}^KM^n_{ij}\ <\ \infty, } 
\end{align*}
where the second step above uses that $S_j(m;n)$ is independent of $\bfZ_0,\ldots,\bfZ_n$ (but not $\bfZ_{n+1},\bfZ_{n+2},\ldots$). Further, the strict inequality 
\begin{equation}
\sum_{n=0}^\infty\sum_{j=1}^K M^n_{ij}<\infty
\end{equation}
follows from $\rho<1$. To prove the necessity, let $\bar{\mu}=\max_1^K\mu_j$. Then by the same reasoning we get 
$\Exp_iT^* \geq \bar{\mu}^{-1} \sum_{n=0}^\infty\sum_{j=1}^K  M^n_{ij} = \infty$ for $\rho \geq 1$ (see Berman and Plemmons \cite{Berman}).   Hence $\Exp_iT^*  < \infty$ is also necessary for $\rho < 1$.   
\end{proof}

\noindent We now have the following corollary to Theorem \ref{thmHA}.
\begin{corollary}
The busy period $T<\infty$ w.p.1 if and only if the matrix $\bfM$ given by
\beqn
M_{ij}=\frac{\lambda_{ij}}{\mu_i},\,\, i,j=1,\ldots, K,
\eeqn
has largest eigenvalue $\rho(\bfM)\leq$ 1.
\end{corollary}
\noindent Similarly, a corollary to Theorem \ref{thmHA2} is stated below.
\begin{corollary}
For the busy period $T$, $\Exp T < \infty$ if and only if $\rho < 1$. 
\end{corollary}

 \subsection{Further applications}
 
Let $B_{i;z}$ denote the length of the busy period initiated by a class $i$ customer with
service requirement $z$ ($B_i$ is that of the standard busy period initiated by a class $i$ customer, that is, taking
$z=S_i$). Let further 
\begin{equation}\label{12.5d}\tau_j =\  \Exp_i\sum_{n=0}^\infty \sum_{m=1}^{Z_n^{(j)}}S_j(m;n)
\end{equation}
be the expected total time in $[0,B_j)$ where the customer being served
is of class $i$. As before, $\bfG$  is  the diagonal matrix with the $\mu_i$ on the diagonal.

\begin{lemma}\label{Lemma22.7a}  \tcr{Assume $\rho<1$.} Then:\\
{\rm (i)} $(\Exp_i\bftau_j)_{i,j=1,\ldots,K}\ =\ (\bfI-\bfM)^{-1}\bfDelta\,; $ {\rm (ii)} $\Exp B_i\ =\ \bfe^\top_i(\bfI-\bfM)^{-1}\bfDelta\bfe\,;$\\ 
{\rm (iii)} $\Exp B_{i;z}\ =\ z\beta_i$ where $\beta_i\, =\,\bfe^\top_i\bfLambda(\bfI-\bfM)^{-1}\bfDelta\bfe\,;$\\
{\rm (iv)} $ B_{i;z}/z\to\beta_i$ in probability as $z\to\infty$.
\end{lemma}
\begin{proof} (i) follows immediately since the $ij$ element of $(\bfI-\bfM)^{-1}\bfDelta$ is
\[\sum_{n=0}^\infty M^n_{ij}/\mu_j\ =\ \Exp_i\sum_{n=0}^\infty Z_n^{(j)}/\mu_j\ =\ \Exp\tau_i,\]
and (ii) follows from (i) by summing over $j$. For (iii) and (iv), we may (by work conservation) assume
that the discipline is preemptive-resume. The workload process during service of a class $i$ customer
evolves as a standard compound Poisson process with arrival rate $\bar{\lambda}^i = \sum_{j=1}^N \lambda_{ij}$
and with cumulative distribution function
\[\sum_{j=1}^K\frac{\lambda_{ij}}{\bar{\lambda}^i}\Prob(B_j\le x)\]
for the jumps. For this system, the rate of arriving work is 
$\bar{\lambda}^i \sum_{j=1}^K\lambda_{ij}/\bar{\lambda}^i\, \Exp B_j$, which is the same as $\beta_i$.
Now we may simply appeal to standard compound Poisson results to obtain (iii) and (iv). This concludes the proof. 
\end{proof}

\section{Busy period results} \label{secBP}\setcounter{equation}{0}
\tcr{In this section, we begin by assuming $\rho(M)\leq 1$.}
Let $B_\bfx$ denote the busy period when the system starts from the state $\bfx \in \mathbb{Z}^K_+$, that is, the time period until the system becomes empty. In particular, when $\bfx$ consists of a single customer of class $i$,
we denote the busy period as $B_i$, and $B_{i,s}$ is the busy period when his remaining service is $s$.
Define $g_\bfx$ to be the LST of $B_\bfx$, i.e.,  
$g_\bfx(\theta) = \Exp_\bfx[\e^{-\theta B_\bfx}]$ for $\bfx \in \mathbb{Z}^K_+$ and similarly for $g_i,g_{i,s}$.

\subsection{The busy period Laplace transform} \label{sec6}


For the $M/G/1$ queue, when $K=1$, it is well known that the LST of the busy period $B$ is given by
\begin{equation}
g(\theta) = \psi(\theta+ \lambda - \lambda g(\theta)),
\end{equation}
where $\psi$ is the LST of the service time and $\lambda$ is the arrival rate.  See, for example, Neuts~\cite{Neuts} or Wolff~\cite{Wolff}.
We shall use the  branching process connection  to derive a similar fixed point equation for our model.

We first observe that the busy period of the system corresponding to an arbitrary initial state $\bfx = (x_1,\dots,x_K) \in \mathbb{Z}^K_+$ is the independent sum of busy periods, each of which corresponds to the branching process starting with a single customer. 
This gives immediately that
\begin{equation}\label{19.1a} g_\bfx(\theta)\ =\ g_1^{x_1}(\theta)\cdots g_K^{x_K}(\theta)\quad\text{when }
\bfx = (x_1,\dots,x_K) \in \mathbb{Z}^K_+
\end{equation}
Hence, it is sufficient to calculate $g_i,g_{i,s}$. Recall
that $\psi_i$ is the LST of the service time distribution $F_i$ of a class $i$ customer.

\begin{theorem} \label{thm-LT} 
For $\theta\ge 0$,
\begin{equation}\label{19.1b} 
g_{i,s}(\theta)\ =\ {\rm exp} \Bigl\{-s\Bigl(\theta+\bar{\lambda}^i-\sum_{j=1}^K\lambda_{ij} g_j(\theta)\Bigr)\Bigr\}.
\end{equation}
Further, 
\begin{equation}\label{19.1c} 
g_i(\theta)\ =\ \psi_i \Bigl(\theta+\bar{\lambda}^i-\sum_{j=1}^K\lambda_{ij} g_j(\theta)\Bigr)\,,\quad i=1,\ldots,K,
\end{equation}
and the vector
$\bigl(  g_1(\theta),\ldots  g_K(\theta)\bigr)$ is the minimal non-negative and non-increasing solution of \tcr{this system of equations.}
\end{theorem}
\begin{proof}
Clearly, $B_{i,s}$ is the service time $s$ plus the busy periods of all customers arriving during service. But the number
of such customers of class $j$ \tcr{is Poisson$(  \lambda_{ij} s)$ and} so their busy periods add up to a compound Poisson random variable
with LST ${\rm exp}\bigl\{\lambda_{ij}  s  \bigl(g_j(\theta)-1\bigr)\bigr\}$. The independence for different $j$ then gives
\[g_{i,s}(\theta)\ =\ \e^{-\theta s}\prod_{j=1}^K {\rm exp}\bigl\{\lambda_{ij}s\bigl(g_j(\theta)-1\bigr)\bigr\},\]
which is the same as \eqref{19.1b}. Integrating with respect to \ $ F_i(\dd s)$ then gives \eqref{19.1c}.

Now consider another non-negative solution $\bigl(\widetilde g_1(\theta),\ldots,\widetilde g_K(\theta)\bigr)$ of \eqref{19.1c}. 
Define the depth $D$
of the multitype Galton-Watson family tree as 
$D\,=\, \max\bigl\{n\ge 0:\,\bfZ_n\ne\bfzero \bigr\}$ and let $g_i^{(n)}(\theta)=\Exp[\e^{-\theta B_i};\,D\le n]$.
Here $D=0$ means no arrivals during service. This occurs with probability $\e^{-\bar{\lambda}^iS_i}$ given $S_i$, and so
$g_i^{(0)}(\theta)=\psi_i(\theta+\bar{\lambda}^i)$. The assumptions on $\widetilde g_j(\theta)$  then gives
$\widetilde g_i(\theta)\ge g_i^{(0)}(\theta)$. Further, the same reasoning as that leading to \eqref{19.1c} gives
\[ g_i^{(n+1)}(\theta)\ =\ \psi_i \Bigl(\theta+\bar{\lambda}^i-\sum_{j=1}^K\lambda_{ij} g_j^{(n)}(\theta)\Bigr)\,.\]
By induction starting from $g_i^{(0)}(\theta)\le \widetilde g_i(\theta)$ we then get $g_i^{(n)}(\theta)\le \widetilde g_i(\theta)$
for all $n$. The proof is completed by observing that $\rho(\bfM)\le 1$ implies $D<\infty$ and hence
$g_i^{(n)}(\theta)\uparrow g_i(\theta)$. 
\end{proof}

\begin{example}
Consider a network with $K=2$ users, $\lambda_{11}=\lambda_{22}=0$ and $F_i$ exponential$(\mu_i)$. 
Then \eqref{19.1c} has the form
\[ g_1\ =\ \frac{\mu_1}{\mu_1+\theta+\lambda_{12}-\lambda_{12}g_2}\,,\quad
 g_2\ =\ \frac{\mu_2}{\mu_2+\theta+\lambda_{21}-\lambda_{21}g_1}\]
where for brevity $g_i$ means $g_i(\theta)$. 
This gives 
\[  g_1\ =\ \frac{
\mu_1 \lambda_{21} - \mu_2 \lambda_{12}  + (\mu_1 + \theta + \lambda_{12} )(\mu_2 + \theta + \lambda_{21} )  - \sqrt{\Delta}
}{2 \lambda_{21} (\mu_1 + \theta + \lambda_{12} ) }, \] 
\[  g_2\ =\ \frac{
 - \mu_1 \lambda_{21} + \mu_2 \lambda_{12} + (\mu_1 + \theta + \lambda_{12} )(\mu_2 + \theta + \lambda_{21} )  - \sqrt{\Delta}
}{2 \lambda_{12} (\mu_2 + \theta + \lambda_{21} ) }, \] 
where
\[
\Delta = [ \mu_1 \mu_2 + \lambda_{12} \lambda_{21} + \theta^2 + \theta(\mu_1 + \mu_2 + \lambda_{12} + \lambda_{21} ) ]^2  - 4 \mu_1 \mu_2 \lambda_{12} \lambda_{21} . 
\]

\end{example}

\subsection{Busy period asymptotics}\label{sec7}

In this section, we offer some observations on the tail asymptotics 
of the busy period in the case of heavy-tailed
service time distributions. For light-tailed service time distributions, we refer the reader to the recent work of Palmowski and  Rolski~\cite{Rolski}.
For the current case of heavy tails,  we refer the reader to Zwart~\cite{Bert}, Jelenkovi\'c and  Momcilovi\'c~\cite{Predrag}
and Denisov and Shneer~\cite{DenisSeva}.

The key idea in both Jelenkovi\'c and Momcilovi\'c \cite{Predrag} and in Zwart \cite{Bert} (as in many other instances of heavy-tailed behavior) is the principle of \emph{one big jump}.
For busy periods, this leads us to expect a large busy period to occur as consequence of one large service time. For concreteness, consider
the standard $M/G/1$ queue with \tcr{$\rho<1$} and suppose there is a single large service time of size $S=z$. The workload
after the large jump is $u+z$ for some small or moderate $u$. The workload then decreases at the rate $1-\rho$
until it reaches 0 and the busy period terminates. By the Law of Large Numbers (LLN), the time of termination is approximately $(z+u)/(1-\rho)$. Since the time before the big jump can be neglected, we have $B>x$ if and only if $z>(1-\rho)x$. Both Asmussen~\cite{SA98} and  Foss and Zachary~\cite{FossZ} show the probability of this large jump is 
\tcr{asymptotically equal to $\Prob\bigl(S>(1-\rho)x\bigr) \Exp\sigma$ for large $x$, where $\sigma$ is  the  number} of customers served in a busy period. 
\tcr{But $\Exp\sigma= \sum_{n=0}^\infty \rho^n = 1/(1-\rho)$. Indeed, 1 corresponds
to the customer initiating the busy period,   $\rho$ is the number of customers arriving while he is in service (the first generation), $\rho^2$
is the number of customers arriving while they are in service,
and so forth. In the framework of branching processes, $\rho^n$ is the number of individuals in the $n$th generation. These considerations} lead to 
\begin{align}
\label{22.7e}  \Prob(B>x)\ &\sim\ \frac{1}{1-\rho} \Prob\bigl(S>(1-\rho)x\bigr),
\end{align}
which Jelenkovi\'c and Momcilovi\'c \cite{Predrag} shows to be the correct asymptotics if the service time distribution is subexponential and
square root insensitive, i.e.\ with a heavier tail than $\e^{-\sqrt{x}}$.

Generalizing this approach to our multiclass system, we recall that 
$\beta_i\, =\,\sum_{j=1}^K\lambda_{ij}\Exp B_j$ and we
introduce a subexponential and
square root insensitive reference distribution  $F$ for which
the individual service time distributions are related as
\begin{equation}\label{0509a}
\Fb_i\bigl(x/(1+\beta_i)\bigr)\ \sim\ c_i\Fb(x).
\end{equation}
In practice, one chooses $\Fb(x)$  as $\sup_i\Fb_i\bigl(x/(1+\beta_i)\bigr)$.
This is common in heavy-tailed studies involving distributions with different degrees of heavy-tailedness. In particular, it allows some $F_j$ to be light-tailed ($c_j=0$). \\
\indent Recalling the interpretation of $\beta_i$ as the rate of arriving work while a class $i$ customer is in service,
a big service time $S_i$ of a class $i$ customer will lead to  $B_i>x$ precisely when $S_i(1+\beta_i)>x$.
Using 
the same reasoning as for~\eqref{22.7e}, we first note that $(\bfM^n)_{ij}$ is the number of type $j$ progeny of a type $i$ ancestor. Hence if $\rho(\bfM) < 1$, the probability that   one of these large service times occur in $[0,B_i)$ is approximately 
\begin{align*}\sum_{n=0}^\infty \sum_{j=1}^K (\bfM^n)_{ij}\Fb_j\bigl(x/(1+\beta_j)\bigr)\ \sim\ d_i\Fb(x),\\ \intertext{where}\
 d_i\,=\,   \sum_{n=0}^\infty \sum_{j=1}^K (\bfM^n)_{ij}c_j 
 \,=\,     \sum_{j=1}^K (\bfI-\bfM)_{ij}^{-1}c_j.  \end{align*}
 \tcrr{Equivalently, the $d_i$ solve
 \begin{align}\label{30.7b} d_i\ &=\ c_i+\sum_{j=1}^K m_{ij}d_j\,.
\end{align}}
As for the standard $M/G/1$ queue, it is straightforward to verify that this is an asymptotic lower bound.
\begin{proposition}\label{Prop:0509a}
Assume that $F$ in \eqref{0509a} is subexponential with finite mean \tcb{so that} $c_k>0$ for some $k$ and \tcr{$\rho(\bfM) < 1$.} Then 
for each $i=1,\ldots,K$,
\begin{equation}\label{0509b}
\liminf_{x\to\infty}\frac{\Prob(B_i>x)}{\Fb(x)}\ \ge\ d_i.
\end{equation}
\end{proposition}

\begin{remark} \rm \tcrr{Square root insensitivity of $F$ is not needed for Proposition~\ref{Prop:0509a}.}
The assumption
\begin{equation}\label{0509avar}
\Fb_i(x) \sim\ \widetilde c_i\Fb_0(x)
\end{equation}
may apriori be more appealing than \eqref{0509a} since it does not involve evaluation of the $\beta_i$.
However, it is closely related. The reason is that if $F$ is regularly varying with $\Fb(x)=L(x)/x^\alpha$, then 
\eqref{0509a} and \eqref{0509avar} with $F_0=F$ are equivalent, with the constants related by
$c_i=\widetilde c_i(1+\beta_i)^\alpha$. For $F_0$ lognormal or Weibull with tail $\e^{-x^\delta}$
(where $\delta<1/2$ in the square root insensitive case), one has, for $\gamma_1>\gamma_2$, $\Fb_0(\gamma_1x)=o\bigl(\Fb_0(\gamma_2x)\bigr)$. Hence if \eqref{0509avar} holds, we may define $\beta^*=\max_1^K\beta_j$
and take $\Fb(x)=\Fb_0\bigl(x/(1+\beta^*)\bigr)$, where $c_j=1$ if $\beta_j=\beta^*$ 
and $c_j=0$ if $\beta_j<\beta^*$.
\end{remark}

\newcommand{\eqdistr}{\stackrel{{\footnotesize \cal D}}{=}}

\indent The $M/G/1$ literature leads to the conjecture that further contributions
to  $\Prob(B_i>x)$ can be neglected, i.e.\ that 
$\Prob(B_i>x)\sim d_i\Fb(x)$  \tcrr{in the square-root insensitive case.}
However, the upper bound \tcrr{is much more difficult (even in the single-class $M/G/1$ setting)}
and follows 
in the regular varying case from more general
results  recently   established in Asmussen \& Foss~\cite{SASF17}:
\begin{theorem}\label{Th:29.7a}
Assume that in addition to the conditions of Proposition~\ref{Prop:0509a} that $F$ is regularly
varying. Then $\Prob(B_i>x)\sim d_i \Fb(x)$
for each $i=1,\ldots,K$. 
\end{theorem}
In the proof, we need:
\begin{lemma}\label{L:30.7a} Let  $S$ be subexponential and let the conditional distribution of $N$
given $S=s$ be \tcb{Poisson$(\lambda s)$}. Then $\Prob(S+N>x)\sim\Prob\bigl(S(1+\lambda)>x\bigr)$ \tcrr{as $x\to\infty$}.
Further, the conditional distribution of $(S,N)/(S+N)$ given $S+N>x$ converges to the one-point distribution
at $((1/1+\lambda),\lambda/(1+\lambda))$.
\end{lemma}
\begin{proof} The argument is standard, with the key intuition being that the variation in $S$ dominates that of 
the Poisson distribution,
so that $N$ can be replaced by its conditional expectation $\lambda S$ given $S$. 
Firstly, note that if $x$ is so large that
$x-x^{1/2}>2\lambda x^{1/2}$ and $N(x^{1/2})$ is Poisson$(\lambda x^{1/2})$, then
\begin{align*}\Prob(S+N>x, S<x^{1/2})\ &\le\ 
\Prob\bigl(x^{1/2}+N(x^{1/2})>x\bigr)\ \le\ \Prob\bigl(N(x^{1/2})>
2\lambda x^{1/2}\bigr),
\end{align*}
which (by large deviations theory) tends to zero  \tcrr{faster than $\e^{-\delta x^{1/2}}$ for some $\delta>0$}, and hence faster than
\tcb{ $\Prob\bigl(S(1+\lambda) > x\bigr)$ }. Secondly, $N/S\to \lambda$ as $y\to\infty$ given $S>y$ and so
\begin{align*}\Prob(S+N>x, S\ge x^{1/2})\ &\sim\ \Prob\bigl(S(1+\lambda)>x,  \tcb{ S\ge x^{1/2} \bigr)},
\end{align*}
the latter equaling $ \Prob\bigl(S(1+\lambda)>x\bigr)$ for large $x$. This proves the first statement,
and the second follows since (asymptotically) only large values of $S$ contribute to large values
of $S+N$, and in this regime $N/S\sim\lambda$. 
\end{proof}

The set-up of ~\cite{SASF17} is a set of random variables $(B_1,\ldots,B_K)$  satisfying
\begin{equation}\label{AF30.6a}B_i \ \eqdistr\  S_i+ \sum_{j=1}^K\sum_{m=1}^{N_{j;i}}B_{m;\tcrr{j}}.\end{equation}
The assumptions for  \eqref{AF30.6a} are that all $B_{m;j}$ are independent of
the vector
$(S_i,N_{1;i},\ldots,N_{K;i})$, that they are mutually independent, and that $B_{m;j}\eqdistr B_j$.
Further, all random variables are non-negative. In our multiclass queue, $B_i$ is the length of the busy period initiated
by a class $i$ customer, $S_i$ is the service time, and $N_{j;i}$ is the number of class $j$ customers
arriving during his service. In the following, we omit the index $i$ and instead express the dependence
on $i$ in terms of a governing probability measure $\Prob_i$.

\begin{proof} \emph{of Theorem~\ref{Th:29.7a}} \ \ 
To apply the results of \cite{SASF17}, we first need to verify a condition on multivariate regular
variation (see~\cite{Elephant} for background) of the vector $(S,N_1,\ldots,N_K)$. Its first part is that 
$\Prob(S+N_1+\cdots+N_K>x)$ $\sim b_i\Fb(x)$ for some $b_i$. This is immediate from 
Lemma~\ref{L:30.7a} by taking $N=N_1+\cdots+N_K$, $\lambda=\overline\lambda_i$, 
$b_i=\widetilde c_i(1+\overline\lambda_i)^\alpha$. A minor extension
of the proof of Lemma~\ref{L:30.7a} further yields that, given $S+N_1+\cdots+N_K>x$,
\bneqn \label{eqangular}
\frac{1}{S+N_1+\cdots+N_K}\bigl(S, N_1, \ldots ,N_K)\ \to 
\frac{1}{1+\overline\lambda_i}\bigl(1,\lambda_{i1}, \ldots , \lambda_{iK}\bigr),
\eneqn
where the limit is taken as $x \rightarrow \infty$.
 This establishes the second part, namely the existence
of the so-called angular measure (in this case a one-point distribution at the right-hand side of (\ref{eqangular})).

It now follows from~\cite{SASF17} that $\Prob(B_i>x)\sim d_i^* \Fb(x)$, where the $d_i^*$ solve
the set of linear equations
\begin{align}\label{30.7a} d_i^*\ &=\ c_i^*+\sum_{j=1}^K m_{ij}d_j,\end{align}
and
\[c_i^*\ =\ \lim_{x\to\infty}\frac{1}{\Fb(x)}\Prob_i(S+N_1\overline r_1+\cdots+N_K\overline r_K>x)\quad
\text{with } \overline r_j=\Exp_j \tcrr{B}.
\] 
Comparing with \eqref{30.7b}, we see that we need only check that $c_i^*=c_i$.
But by similar arguments to those above, 
\begin{align*}
&\Prob_i(S+N_1\overline r_1+\cdots+N_K\overline r_K>x)\ \sim\
\tcb{ \Prob(S(1+\lambda_{i1}\overline r_1+\cdots+\lambda_{iK} \overline r_K ) >x) } \\ &=\ 
\Prob(S(1+\beta_i)>x)\ \sim\ \widetilde c_i(1+\beta_i)^\alpha\Fb(x)\ =\ c_i\Fb(x),\end{align*}
where \tcrr{parts (ii) and (iii) of} Lemma~\ref{Lemma22.7a} are employed in the second step. 
\end{proof}

\begin{remark} \rm
The general subexponential case seems much more difficult. One obstacle is that
theory and applications of multivariate subexponentiality is much less developed than for 
the regular varying case. See, however, Samorodnitsky and Sun~\cite{Genna} for a recent
contribution and for further references.
\end{remark}

\section{Conclusion} \label{sec8}\setcounter{equation}{0}
We have introduced a multiclass single-server queueing server in which the arrival rates depend on the current job in service. The model departs from existing state-dependent models in the literature in which the parameters depend primarily on the number of jobs in the system rather than the job in service. \\
\indent The main contributions of this paper can be summarized as follows. Firstly, we formulate the multiclass queueing model and its corresponding fluid model, and provide motivation for its practical importance.  The necessary and sufficient conditions for stability of the queueing system are obtained via the corresponding fluid model. Secondly, by appealing to the natural connection with multitype Galton-Watson processes, we utilize Laplace-Stieltjes transforms to characterize the busy period of the queueing system. Thirdly, we 
present \tcrr{ a preliminary study of} busy period tail asymptotics for heavy-tailed service time distributions
\tcrr{ and give  a complete set of results for the regularly varying case,
using recent results of Asmussen \& Foss~\cite{SASF17}}. Tail asymptotics in our multiclass setting for non-regularly varying heavy-tailed service time distributions, as well as for light-tailed service time distributions, are much more difficult and will be attempted in a separate manuscript.\\

\noindent \textbf{Acknowledgments} 
We are very grateful to a referee for pointing out a problem in our initial proof of the upper bound in Section~\ref{sec7}. We also thank a second referee \tcrr{and an associate editor} for many useful suggestions. 
The first author thanks Dr. Quan Zhou and Professor Guodong Pang for helpful conversations.


\end{document}